\documentclass[11pt]{article}
\usepackage[dvips]{graphicx}
\usepackage{latexsym}
\usepackage{amssymb}

\newcommand{\C}{\mathbb{C}}
\newcommand{\CP}{\mathbb{CP}}

\newcommand{\R}{\mathbb{R}}

\renewcommand{\d}{\mathrm{d}}

\newcommand{\koniec}{\begin{flushright}  $\Box $ \end{flushright}}
\def\be{\begin{equation}}
\def\ee{\end{equation}}

\def\Sm{\Sigma}

\def\r{\rho}

\def\Om{\Omega}

\def\O{\cal O}
\def\om{\omega}

\def\ov{\overline}

\def\wt{\overline{w}}
\def\zt{\overline{z}}

\def\tw{\tilde w}
\def\tz{\tilde z}
\def\lt{\tilde{\lambda}}
\def\p{\partial}
\def\pt{\tilde \p}

\def\ov{\overline}

\def\Kt{\tilde{K}}

\def\a{\alpha}
\def\l{\lambda}

\def\O{{\cal O}}

\def\dom{\delta \Om}

\def\Kt{\widetilde{K}}
\newtheorem{theo}{Theorem}[section] 
\newtheorem{prop}[theo]{Proposition}  
\newtheorem{lemma}[theo]{Lemma}

\newtheorem{col}[theo]{Corollary}

\begin{document}

\title{Einstein--Weyl structures from Hyper--K\"ahler metrics with conformal Killing vectors}
\author{Maciej  Dunajski\thanks{email: dunajski@maths.ox.ac.uk}\\
Paul Tod\\ The Mathematical Institute,
24-29 St Giles, Oxford OX1 3LB, UK}  
\maketitle
\noindent
\abstract {We consider four (real or complex) dimensional
hyper-K{\"a}hler metrics with a conformal symmetry $K$. The three-dimensional 
space of orbits of $K$ is shown to have an Einstein--Weyl structure
which admits 
a shear-free geodesics congruence for which the twist is a constant multiple
of the  divergence. In this case
the Einstein--Weyl equations reduce down to a single second order PDE
for one function. The Lax representation, Lie
point symmetries,  hidden  symmetries and the recursion
operator associated with this PDE are found, 
and some group invariant solutions are considered.}

\section{Three-dimensional Einstein--Weyl spaces}

Three-dimensional Einstein--Weyl (EW) geometries were first considered by Cartan \cite{C43}
and then rediscovered by Hitchin \cite{H82} in the context of twistor theory.
They constitute an interesting generalisation of (the 
otherwise locally trivial) Einstein
condition in three dimensions.

In this paper  we shall consider four-dimensional anti-self-dual (ASD)
vacuum (or complexified hyper-K\"ahler) spaces with a conformal symmetry. 
By a general construction \cite{JT85}
such spaces will give rise to Einstein--Weyl structures on the space
of trajectories of the given conformal symmetry $K$. The cases where
$K$ is a pure Killing vector 
or a tri-holomorphic homothety have been
extensively studied \cite{BF82,W90,CTV96,L91}. Therefore we shall
consider the most general case of $K$ being a conformal,
non-triholomorphic Killing vector. 
We begin by collecting various definitions and formulae concerning
 three-dimensional
Einstein--Weyl spaces (see \cite{PT93} for a fuller account).
In the next section  
we shall give the canonical form of an allowed conformal Killing
vector in a natural coordinate system associated with the K{\"a}hler
potential. 
Then we shall look at solutions to a non-linear Monge--Ampere equation
(the so called `first heavenly equation' \cite{Pl75})
(\ref{pleb1}) for the K\"ahler potential  
which admit the symmetry $K$. This will give rise to
a new integrable system in three dimensions and to the corresponding
EW geometries. In Section \ref{Lax_representation} we shall give the
Lax representation
of the reduced equations. When Euclidean reality conditions are
imposed (Section \ref{special_cases})
we shall recover some known results \cite{BF82,W90} as limiting
cases of our construction.
In Section \ref{lie_symmetries} we shall find 
and classify the Lie point symmetries
of  the field equations in three
dimensions (and so the Killing vectors of the associated 
Weyl structure), and consider  some group invariant solutions.
In Section \ref{hidden_symmetries} we shall study hidden 
symmetries and the recursion
operator associated to the three-dimensional system. 
In Section 8 we shall show that the EW structures studied in this paper
admit a shear free geodesic congruence for which twist
and divergence are linearly dependent.

Let ${\cal W}$ be an $n$-dimensional complex manifold,
with a torsion-free connection
$D$ and a conformal metric $[h]$.
We shall call  ${\cal W}$ a Weyl space if
the null geodesics of $[h]$ are also geodesics for $D$. 
This condition is equivalent to 
\be
\label{ew1}
D_ih_{jk}=\nu_ih_{jk}
\ee
for some one form $\nu$. Here $h_{jk}$ is a representative metric in
the conformal class. The indices 
$i, j, k, ...$ run from 1 to $n$.
If we change this representative by
$h\longrightarrow \phi^2 h$, then $\nu\longrightarrow
\nu+2\d\ln{\phi}$.  
The one form $\nu$ `measures' the difference between $D$ and
the Levi-Civita connection $\nabla$ of $h$:
\[
(D_i-\nabla_i)V^j=\gamma^j_{ik}V^k,\qquad\mbox{where}\qquad
\gamma^i_{jk}:=-\delta_{(j}^i\nu^{}_{k)}+\frac{1}{2}h_{jk}h^{im}\nu_m.
\]
The Ricci tensor $W_{ij}$ of $D$ is related to the Ricci tensor
$R_{ij}$ of $\nabla$ by
\[
W_{ij}=R_{ij}+\frac{n-1}{2}\nabla_i\nu_j-\frac{1}{2}\nabla_j\nu_i
+\frac{n-2}{4}\nu_i\nu_j+h_{ij}\Big(-\frac{n-2}{4}\nu_k\nu^k+\frac{1}{2}
\nabla_k\nu^k\Big).
\]
The relation between the curvature scalars is
\[
W:=h^{ij}W_{ij}=R+(n-1)\nabla^k\nu_k-\frac{(n-2)(n-1)}{4}\nu^k\nu_k.
\]
The conformally invariant Einstein--Weyl (EW) condition on
$({\cal W}, h, \nu)$ is
\[
W_{(ij)}=\frac{1}{n}W h_{ij}.
\]

From now on we shall assume
that $dim{\cal W}=3$. The Einstein--Weyl equations can then be 
written
\be
\label{ew2}
\chi_{ij}:=R_{ij}+\frac{1}{2}\nabla_{(i}\nu_{j)}+\frac{1}{4}\nu_i\nu_j
-\frac{1}{3}\Big(R+\frac{1}{2}\nabla^k\nu_k+\frac{1}{4}\nu^k\nu_k\Big)h_{ij}=0.
\ee
Here $\chi_{ij}$ is the trace-free part
of the Ricci tensor of the Weyl connection.
In three dimensions  the general solution of (\ref{ew1})-(\ref{ew2}) 
depends on four arbitrary functions
of two variables \cite{C43}.  
In this paper we shall consider a class of solutions
to the EW equations which depend on 
two arbitrary functions
of two variables.

All three-dimensional EW spaces can be obtained as spaces of trajectories of
conformal Killing vectors in four-dimensional ASD manifolds:
\begin{prop}[Jones \& Tod \cite{JT85}]
\label{prop_JT}
Let $({\cal M}, g)$ be an ASD four manifold with a conformal Killing
vector $K$.
An EW structure on the space  ${\cal W}$ of trajectories of $K$ 
{\em (}which is assumed to be non-pathological{\em)} is
defined by
\be
\label{EWs}
h:=|K|^{-2}g-|K|^{-4}{\bf K}\odot {\bf K},\;\;\; \nu:=2|K|^{-2}
\ast_g({\bf K}\wedge \d{\bf K}),
\ee
where $|K|^2:=g_{ab}K^aK^{b}$, ${\bf K}$ is the one form dual to $K$ and
$\ast_g$ is taken with respect to $g$. All EW structures arise in this way.
Conversely, let $(h, \nu)$ be a three--dimensional EW structure 
on ${\cal W}$, and
let $(V, \omega)$ be a function and a one-form on ${\cal W}$ which satisfy
the generalised monopole equation
\be
\label{EWmonopole}
\ast_h(\d V+(1/2)\nu V) =\d\omega,
\ee
where $\ast_h$ is taken with respect to $h$. Then
\be
g=V^2h+(\d t+\omega)^2
\ee
is an ASD metric with an isometry $K=\p_t$.
\end{prop}
\def\wt{\overline{w}}
\def\zt{\overline{z}}

\section{Hyper-K\"ahler metrics with conformal Killing vectors}
\label{HS_conformal}
Let $g$ be a compexified  hyper-K\"ahler (i.e. ASD vacuum)
metric on a complex
four-manifold ${\cal  M}$ and  $x^{AA'}=(w, z, \tw, \tz)$ be a 
null coordinate system on ${\cal  M}$. Locally $g$ can be put in the form 
\be
\label{pmetric}
\d s^2=\Om _{w\tw}\d w \d \tw+\Om _{w\tz}\d w \d \tz
+\Om _{z\tw}\d z \d \tw +\Om _{z\tz}\d z \d \tz
\ee
(subscripts denote partial differentiation) where $\Om=\Om(w, z, \tw, \tz)$ is solution of the first heavenly
equation \cite{Pl75}
\be
\label{pleb1}
\Om_{w\tz}\Om_{z\tw}-\Om _{w\tw}\Om_{z\tz}=1.
\ee 
Assume that $g$ admits  a conformal Killing vector $K$; i.e. 
${\cal L}_Kg=\eta g$, or equivalently
\[
\nabla_aK_b=\phi_{A'B'}\varepsilon_{AB}+\psi_{AB}\varepsilon_{A'B'}
+(1/2)\varepsilon_{A'B'}\varepsilon_{AB}\eta
\]
where symmetric spinors 
$\phi_{A'B'}$ and $\psi_{AB}$ are respectively self-dual and 
anti-self-dual parts of 
the covariant derivative of $K$.
The well known formula $\nabla_a\nabla_bK_c=R_{abcd}K^d$ relating the
second covariant derivative of $K$ to the Riemannian curvature implies
that in vacuum
\[
\nabla_{AA'}\phi_{B'C'}=2C_{A'B'C'D'}{K_A}^{D'}-2\varepsilon_{A'(B'}
\nabla_{C')A}\eta, \;\;\nabla_{AA'}\nabla_{BB'}\eta=0,
\;\; C_{ABCD}{\nabla^{A}}_{A'}\eta=0.
\]
Here $C_{ABCD}$ and $C_{A'B'C'D'}$ are respectively the ASD and SD Weyl
spinors. In particular in an ASD vacuum $\phi_{A'B'}=const$ and $\eta=const$
(or the space time is of type $N$). In this paper we shall analyse a
situation where $K$ is not hyper--surface orthogonal and
$det(\phi_{A'B'})\neq 0, \eta\neq 0$.
\begin{lemma}
\label{K_lemma}
In an ASD vacuum the most general conformal 
Killing vector with det$(\phi_{A'B'})\neq
0$ can be transformed to the form
\be
\label{can_Killing}
K=\eta (z\p_z+\tz\p_{\tz})+\rho (z\p_z-\tz\p_{\tz}).
\ee
\end{lemma}
{\bf Proof.}
In the adopted coordinate system a basis of SD two-forms is
\[
\Sm^{0'0'}=\d \tw\wedge \d \tz,\;\;\Sm^{1'1'}=\d w\wedge \d z,
\]
\[
\Sm^{1'0'}=\Om _{w\tw}\d w\wedge\d \tw+\Om _{w\tz}\d w\wedge\d \tz+
\Om _{z\tw}\d z \wedge\d \tw+\Om _{z\tz}\d z\wedge\d \tz.
\]
Let $K=K^A\p/\p w^A+{\Kt}^A\p/\p \tw^A$, where $w^A=(w, z)$ and
$\tw^A=(\tw, \tz)$. The action of $K$ on self-dual two forms
is determined by
\begin{eqnarray*}
{\cal L}_K\Sm^{0'0'}&=& m\Sm^{0'0'}+n\Sm^{1'1'},\\
{\cal L}_K\Sm^{1'1'}&=& \tilde{n}\Sm^{0'0'}+\tilde{m}\Sm^{1'1'},\\
{\cal L}_K\Sm^{0'1'}&=& \eta \Sm^{0'1'}
\end{eqnarray*}
for some constants $m, \tilde{m}, n, \tilde{n}$. 
This is because for non-degenerate $\phi_{A'B'}$ the K{\"a}hler
structure can be identified with $\d K_+=\phi_{A'B'}\Sm^{A'B'}$.
It follows that 
$n=\tilde{n}=0$, and $K^A=m w^A,\;\Kt^A =\tilde{m}\tw^A$. 
From $2\Sm^{0'0'}\wedge\Sm^{1'1'}=\Sm^{0'1'}\wedge\Sm^{0'1'}$ 
we find  that $\eta:=(m+\tilde{m})/2$. Define $\rho:=(m-\tilde{m})/2$.
We
have the  freedom to transform $w^A \rightarrow {W}^A(w^B)$ and 
$\tw^A \rightarrow \widetilde{W}^A(\tw^B)$ in a way which preserves
$\Sm^{0'0'}$ and $\Sm^{1'1'}$. Put
$Z=z^2/2, W=w/z, \widetilde{Z}=\tz^2/2, \widetilde{W}=\tw/\tz$. 
This yields (coming back to $(w^A, \tw^A)$) (\ref{can_Killing}).
Now
\[
\nabla_{AA'}{K^A}_{B'}=\left (
\begin{array}{cc}
0&\rho+\eta\\
\rho-\eta &0
\end{array}
\right ). 
\]
\koniec
The real form of the 
Killing vector (\ref{can_Killing}) also appears in the list of Lie
point symmetries of (\ref{pleb1}) given in \cite{BW89}.
\subsection{Symmetry reduction}
In this section we shall look at the heavenly equation (\ref{pleb1})
with the additional constraint ${\cal L}_Kg=\eta g$. This will lead to a
new integrable equation describing a class of three--dimensional
Einstein--Weyl geometries.
\begin{prop}
Every ASD vacuum metric
with conformal symmetry is locally given by
\begin{eqnarray}
\label{EWmetric}
\d s^2&=&e^{\eta t}(V^{-1}h+V(\d t+\om)^2)\;\;\;\;\;\mbox{{\em where}}\\
 h&=&-e^{2\rho u}\d w \d \tw-\frac{1}{16}(\eta^2F\d u+\eta(F_w\d
 w-F_{\tw}\d \tw)-\d F_u)^2\\
\om &=&\frac{(\eta F_w-F_{uw})\d w+(\eta F_{\tw}+F_{u\tw})\d \tw}
{\eta^2F-F_{uu}},\;\;\;\;\;V=\frac{1}{4}(\eta^2F-F_{uu}),
\end{eqnarray}
and $F=F(w, \tw, u)$ is a holomorphic function on an open set
${\cal W}\subset \C^3$ which satisfies  
\be
\label{EWred}
(\eta F_{\tw}+F_{u\tw})(\eta F_w-F_{uw})-(\eta^2F-F_{uu})F_{w\tw}
=4e^{2\rho u}
\ee
for constants $\eta, \rho\in  \C$.
\end{prop}
\begin{col}
The metric $h$ is defined on the space ${\cal W}$
of trajectories of $K$ in ${\cal M}$.
From Proposition {\em \ref{prop_JT} } it follows 
that $h$ is the most general
EW
metric which arises as a reduction of ASD vacuum solutions
by a conformal Killing vector. 
Equation {\em (\ref{EWred})} is therefore equivalent to the Einstein--Weyl
equations {\em (\ref{ew2})}.
\end{col}
{\bf Proof.}
The general ASDV metric can  locally be given by (\ref{pmetric}).
From Lemma \ref{K_lemma} it follows that we can take $K$ as in (\ref{can_Killing}).
Perform the coordinate transformation $(z, \tz)\rightarrow(t, u)$
given by 
\[
2t:=\ln(z^{1/m}\tz^{1/{\tilde{m}}}),\;\;
\;2u:=\ln(z^{1/m}\tz^{-1/{\tilde{m}}}).
\]
In these coordinates $K=\p_t$ and  so
$\Om(t, u, w, \tw)=e^{\eta t}F(u, w, \tw)$. The first heavenly
equation is equivalent to (\ref{EWred}). Rewriting the metric (\ref{pmetric}) 
in the new coordinate system yields (\ref{EWmetric}) and 
$det(h)=-(1/4)V^2e^{4\rho u}$.

The dual to $K$ is  ${\bf K}=e^{\eta t}{V}(\d t+\om)$.
From Proposition \ref{prop_JT} we find the  EW one-form  to be
\begin{eqnarray*}
\nu&=&2\ast_g\frac{{\bf K}\wedge \d {\bf K}}{|K|^2}=2e^{\eta t}V\ast_g((\d
t +\om)\wedge \d \om)\\
&=&4\rho\d u+
\frac{(2\eta+4\rho)(\eta F_w-F_{uw})\d w+
(2\eta-4\rho)(\eta F_{\tw}+F_{u\tw})\d \tw }{\eta^2F-F_{uu}}
\end{eqnarray*}
where $\ast_g$ is the Hodge operator determined by $g$.\koniec
\section{Lax representation}
\label{Lax_representation}
In this section we shall represent equation (\ref{EWred})
as the integrability condition for a linear system of equations.
We shall interpret the Lax pair  as a  (minitwistor) distribution on 
a reduced projective spin bundle. 
The Lax pair for the first heavenly
equation 
\begin{eqnarray}
\label{LaxH}
L_0:&=&\Om_{w\tw}\p_{\tz}-\Om_{w\tz}\p_{\tw}-\l\p_w,\nonumber\\
L_1:&=&\Om_{z\tw}\p_{\tz}-\Om_{z\tz}\p_{\tw}-\l\p_z
\end{eqnarray}
is defined on the five complex dimensional correspondence space 
${\cal F}={\cal M}\times\CP^1$. Here $\l\in \CP^1$ parametrises null 
self-dual surfaces passing through a point in $\cal M$.
Equations $L_0\Psi=L_1\Psi=0$ have solutions in ${\cal F}$ provided that
$\Om$ satisfies the first heavenly equation (\ref{pleb1}).
The formulation (\ref{LaxH}) is crucial to the twistor construction, as
the projective twistor space on ${\cal M}$ arises as a factor space
of ${\cal F}$ by the distribution $\{L_0, L_1\}$.

Let $\pi_{A'}=(\pi_{0'}, \pi_{1'})$ be coordinates on the fibers of a
bundle $S_{A'}$ of primed spinors. The space ${\cal F}$ can be
regarded as the projectivised version of $S_{A'}$ in a sense that
$\l=\pi_{0'}/\pi_{1'}$.
Define the Lie lift of a Killing vector  $K$ to   
${\cal F}$ by
\be
\label{Klift}
\Kt:=K+Q\p_{\l},\;\;\;\;\mbox{where}\;\;\;Q:=
\pi_{A'}\pi_{B'}\phi^{A'B'}/(\pi_{1'})^2.
\ee
The flow of $\Kt$ in ${\cal F}$ determines the behaviour of
$\a$-planes under the action of $K$ in ${\cal M}$. 
The linear system $L_A$ for equation (\ref{pleb1}) is given by 
(\ref{LaxH}).
The vector fields $(L_0, L_1, \Kt)$ span an integrable 
distribution. This can be seen as follows: 
\begin{eqnarray*}
[K, L_A]&=&-\pi^{A'} ({ {\phi}_{A'} }^{B'}{\epsilon_A}^B+
{{\psi}_{A}}^{B}{\epsilon_{A'}}^{B'}+(1/2)\eta{\epsilon_A}^B
{\epsilon_{A'}}^{B'})\nabla_{BB'}\\
&=&-\pi^{C'}{{\phi}_{C'}}^{A'}\nabla_{AA'}+
({{\psi}_{A}}^{B}+\eta{\epsilon_A}^B)L_B.
\end{eqnarray*}
The Lie  lift of $K$ to $S^{A'}$ is 
\be
\Kt=K+\pi_{A'}\phi^{A'B'}\frac{\p}{\p \pi^{B'}}+
\frac{1}{2}\eta\pi^{A'}\frac{\p}{\p\pi^{A'}},
\ee
so that $[\Kt, L_A]=0$ modulo $L_A$. 

 The projection of $\Kt$ to ${\cal F}$
is given by (\ref{Klift}), where the factor $\pi_{1'}^2$ is used to
dehomogenise a section of $\O(2)$.
If $K$ is given by (\ref{can_Killing}) then $\Kt=K+\rho\l\p_{\l}$.
Introduce an invariant spectral parameter $\lt$ 
(which is constant along $\Kt$) by
$(\l, t)\longrightarrow (\lt:=\l e^{-\rho t}, \hat{t}:=t)$. In the new
coordinates
\[
\p_t=\p_{\hat{t}} -\r\lt\p_{\lt},\;\; 
\p_{\l}=e^{-\r\hat{t}}\p_{\lt},\;\;\;\;\;\; \mbox{so
that}\;\;\Kt=\p_{\hat{t}}.
\]
The linear system for the reduced equation is obtained from 
(\ref{LaxH}) by rewriting it in $(w, \tw, u, \hat{t}, \lt)$ 
coordinates and ignoring $\p_{\hat{t}}$.
This yields (after rescaling)
\begin{eqnarray}
\label{Laxred}
L_{0'}&=&me^ {\tilde{m}u}\Big(F_{w\tw}\Big(\frac{\p}{\p u}+\r
\lt\frac{\p}{\p \lt}\Big)+(\eta F_{w}-F_{uw})\frac{\p}{\p \tw}\Big)+2\lt 
\frac{\p}{\p w}\\
L_{1'}&=&\tilde{m}e^{\tilde{m}u}\Big( (\eta F_{\tw}+ F_{u\tw})
\Big(\frac{\p}{\p u}+\r \lt\frac{\p}{\p \lt}\Big)+(\eta^2 F-F_{uu})\frac{\p}{\p
\tw}\Big)+2\lt \Big(\frac{\p}{\p u}-\r \lt\frac{\p}{\p \lt}\Big)\nonumber.
\end{eqnarray}
The mini-twistor space corresponding to solutions
of (\ref{EWred}) is the quotient of ${\cal F}$ by the integrable
distribution $(L_{0'}, L_{1'}, \Kt)$. 

The existence of a minitwistor
distribution follows from Hitchin's construction \cite{H82};
the basic mini-twistor correspondence
states that points in ${\cal W}$ correspond in ${\cal Z}$
to rational curves with normal bundle ${\cal O}(2)$.  
Let $l_x$ be the line in ${\cal Z}$ that
corresponds to $x\in {\cal W}$. The normal bundle to $l_x$ consists of
tangent vectors  at $x$ (horizontally lifted to $T_{(x,\lambda)}{\cal
  F}_{W}$) modulo the twistor distribution. Therefore we have a sequence
of sheaves over $\C P^1$
\[
0\longrightarrow D_W \longrightarrow \C^3 \longrightarrow
{\cal O}(2)\longrightarrow 0.
\]
We shall identify $T^i{\cal W}\approx S^{(A'}\otimes S^{B')}$.
The map $\C^3 \longrightarrow {\cal O}(2)$ is given by
$V^{A'B'}\longrightarrow V^{A'B'}\pi_{A'}\pi_{B'}$. 
 Its kernel consists of
vectors of the form $\pi^{(A'}v^{B')}$ with $v^{B'}$ varying. The
twistor distribution is therefore $D_W=O(-1)\otimes S^{A'}$ and so $L_{A'}$
is the global section of $\Gamma(D_W\otimes {\cal O}(1)\otimes
S_{A'})$.
Let $Z$ be a totally geodesic two-plane corresponding to a point $Z$ of
a mini-twistor space. This two plane is spanned by vectors of the form
$V^a=\pi^{(A'}v^{B')}$ with $\pi^{A'}$ fixed.
Let $W^a= \pi^{(A'}w^{B')}$ be another vector tangent to  $Z$. The
Frobenius theorem implies that the Lie bracket $[V, W]$ 
must be tangent to some geodesic in $Z$, i.e.
$
[V, W]=aV+bW
$
for some $a, b$. The last equation determines the mini-twistor
distribution $L_{A'}$ to be a horizontal lift of $\pi^{B'}D_{A'B'}$
to the weighted spin bundle by demanding $L_{A'}\pi_{C'}=0$.
The integrability conditions imply
$[L_{A'}, L_{B'}]= 0$, (mod $L_{A'}$).
In fact if one picks two independent 
solutions of a `neutrino' equation on the EW
background, say $\rho^{A'}$ and  $\lambda^{A'}$, then
$\widehat{L}_{0'}:=\rho^{A'}L_{A'}$, and 
$\widehat{L}_{1'}:=\l^{A'}L_{A'}$ commute exactly:
$
[\widehat{L}_{0'}, \widehat{L}_{1'}]=0.
$

\section{Reality conditions}
\label{reality_conditions}
To obtain real Einstein--Weyl metrics we have to impose reality conditions
on the coordinates $(w, \tw, z, \tz)$
\begin{itemize}
 \item The reduction from the Euclidean slice  $(\zt=\tz, \wt=-\tw)$ 
yields positive definite EW metrics with $u:=iv$ for $v\in \R$. 
Without loss of generality we can
impose the condition $m\tilde{m}=1$, so $\eta=\cos\a$ and $\r=i\sin\a$.
The Euclidean version of (\ref{EWred}) is then
\be 
\label{euclid}
(F\cos^2\a +F_{vv})F_{w\wt}-( F_{\wt}\cos\a-iF_{v\wt})( F_w\cos\a+iF_{vw})
=4e^{-2v\sin\a}.
\ee
To obtain another  form  
introduce $G$ by $G=e^{v\sin\a }F$. The transformed equation, the
metric (rescaled by $e^{2v\sin\a}$) and the EW one-form are:
\be
\label{alterform}
(G+G_{vv}-2 G_v\sin\a)G_{w\wt}-
(e^{i\a}G_{\wt}-iG_{v\wt})(e^{-i\a}G_{w}+iG_{vw})=4,
\ee
\begin{eqnarray}
\label{altermetric}
h&=&\d w\d \wt+\frac{1}{16}(G\d v+\d G_v -2G_v\sin\a \d v-ie^{-i\a}G_w\d w
+ie^{i\a}G_{\wt}\d \wt)^2\nonumber\\
\nu&=&-2\sin{\a}\d v+ 
\frac{(2+2\sin^2{\a})(G_{w}\d w+G_{\wt}\d \wt) +i\sin{2\a}(G_{w}\d w-
G_{\wt}\d \wt)}{G+G_{vv}-2G_v\sin\a}\\
&-&\frac{2(\cos{\a}+2i\sin{\a})G_{vw}\d w+ 2(\cos{\a}-2i\sin{\a})G_{v\wt}\d
\wt }{G+G_{vv}-2G_v\sin\a}.
\end{eqnarray}
\item On an ultra-hyperbolic slice we have 
$\zt=\tz, \wt=\tw$ which again implies $u=iv$. The metric (\ref{EWmetric}) 
has signature $(++-)$. Another possibility
is to take all coordinates as real. This gives a different real metric
of  signature $(++-)$. The function $F$ is
real and $\eta=\sinh\a, \r=\cosh\a$.
\end{itemize}
The analogous reality conditions are imposed on the linear system
(\ref{Laxred}).
From now on we shall be mostly concerned with the positive definite case.
The correspondence space is now viewed as a real six-dimensional
manifold. The real lift of a Killing vector is
$
\Kt=\p_t+i\sin\a (\l\p_{\l}-\ov{\l}\p_{\ov{\l}}).
$
\section{Special cases}
\label{special_cases}
Solutions to (\ref{euclid}) describe the most general $EW$ metrics which arise
as reductions of hyper-K\"ahler structures. In this section we look 
at limiting cases and recover hyper-$CR$ EW spaces \cite{CTV96}, and
LeBrun-Ward  EW spaces which come from the $SU(\infty)$ Toda equation.
The real form of the Killing vector (\ref{can_Killing}) is a linear combination of a rotation
and a dilation; 
\[
K=K_D\cos\a+K_R\sin\a,\qquad\a\in[-\pi/2, 0],\;
K_D:=z\p_z+\zt\p_{\zt},\;K_R:=i(z\p_z-\zt\p_{\zt}).
\]
\subsection{LeBrun--Ward spaces}
Take ${\a=-\pi/2}$. Then $K$ is a pure Killing vector which does not
preserve the complex structures on ${\cal M}$. This case was
studied in \cite{BF82,W90,L91}.
Put $F_v=j, F_{\wt}=\ov{p}$ and rewrite equation (\ref{euclid}) as
\begin{eqnarray}
\label{BFhod}
\d \ov{p} \wedge \d j \wedge \d \wt &=&4e^{2v}\d w\wedge\d\wt\wedge \d v\nonumber\\
\d j\wedge \d w\wedge \d v&=&\d \ov{p}\wedge \d \wt \wedge \d w.
\end{eqnarray}
Use $(j, w, \wt)$ as coordinates and eliminate $\ov{p}$ to obtain
\be
\label{Boyer_Finley}
v_{w\wt}+2(e^{2v})_{jj}=0 
\ee
which is the $SU(\infty)$ Toda (or Boyer--Finley) equation \cite{BF82}.
The metric (\ref{EWmetric}) reduces to
\[
h=e^{2v}\d w\d \wt+\frac{1}{16}\d j^2,\;\;\;\;\;\; \nu=4v_j\d j.
\]
This class of EW spaces is characterised by the existence of a 
twist-free, shear-free geodesic congruence \cite{T95}.

Let us come back to complex coordinates and
put $w=e^{s+ \theta}, \tw=e^{s- \theta}$ and $M=2u+2s$.
In the $(s, u, \theta)$ 
coordinates equation (\ref{Boyer_Finley}) and the metric  become
\[
M_{ss}-M_{\theta \theta}-8(e^M)_{jj}=0,\;\;\;\;\; h=-e^{M}(\d s^2- \d
\theta^2)-\frac{1}{16}\d j^2.
\]
Imposing a symmetry in $\theta=\ln(\sqrt{w/\tw})$ direction we arrive at
\[
M_{ss}-8(e^M)_{jj}=0,
\]
which was solved by Ward \cite{W90} who transformed it to a linear
equation. The conclusion is that LeBrun--Ward EW metrics with
$w\p_w-\tw\p_{\tw}$ symmetry are solved by the same ansatz as those
with $\p_w-\p_{\tw}$ symmetry. In Subsection 
\ref{group_solutions} it will be shown
that imposing $\p_w-\p_{\tw}$ symmetry leads to a linear equation even if $\a$
is arbitrary.
\subsection{Hyper-CR spaces} Put ${\a=0}$.
Then $K$ is a triholomorphic conformal symmetry.
The corresponding EW metrics were in \cite{GT98} called `special',
and then referred to as hyper CR (since each complex 
structure on ${\cal M}$  defines a CR structure on ${\cal W}$).
They are characterised by the existence of a sphere of shear-free and
divergence-free geodesic congruences.
The equation (\ref{euclid}) reduces to
\be
\label{conpl}
F_{w\wt}(F+F_{vv})-(F_w+iF_{vw})(F_{\wt}-iF_{v\wt})=4,
\ee
which is the form given in \cite{TD97}. The corresponding Lax pair is
\begin{eqnarray*}
L_{0'}&=&e^ {iv}\Big(iF_{w\wt}\frac{\p}{\p v}-(
F_{w}+iF_{uw})\frac{\p}{\p \wt}
\Big)+2\l 
\frac{\p}{\p w}\\
L_{1'}&=&e^{iv}\Big( (F_{v\wt}+i F_{\wt})
\frac{\p}{\p v}-(F+F_{vv})\frac{\p}{\p
\wt}\Big)-2i\l \frac{\p}{\p v}.
\end{eqnarray*}
\section{Lie point symmetries}
\label{lie_symmetries}
In order to find the Lie algebra of infinitesimal symmetries of
(\ref{alterform}) we shall convert 
it to system of differential forms. 
Introduce $Q$ and $J$ by $J:=G_{\wt},\;\;Q:=(e^{i\a} G-iG_v)$
\begin{eqnarray}
\label{ideal}
\om_1&:=&i\d Q\wedge\d J\wedge \d \wt+e^{-i\a}(Q\d J-J\d Q)\wedge \d
\wt \wedge \d v\nonumber\\
& &+\d Q\wedge\d \ov{Q}\wedge \d v-4\d w\wedge \d
\wt\wedge \d v,\nonumber\\
\om_2&:=&\d Q\wedge\d w\wedge\d v +e^{i\a}J\d w\wedge \d \wt\wedge \d
v-i\d J\wedge\d w\wedge \d \wt.
\end{eqnarray}
This system forms a closed differential ideal. Its integral manifold 
is a subspace of $\R^6$ on which $\om_{\mu}=0$. This integral
manifold represents a solution to (\ref{alterform}).

 Let $X$ be a vector field on  $\R^6$. The action of $X$ does not
change the integral manifold if $
{\cal L}_X\om_{\mu}=\Lambda_{\mu}^{\nu}\om_{\nu}
$
where $\mu, \nu =1, 2$ and  
$\Lambda_{\mu}^{\nu}$ is a matrix of differential forms. 
The general solution is
\begin{eqnarray*}
X&=&(Aw+B)\frac{\p}{\p w}+(\ov{A}\wt+\ov{B})\frac{\p}{\p {\wt}}
+C\frac{\p}{\p v}+\frac{1}{2}(A+\ov{A}) G\frac{\p}{\p G}\\
&+&D_1e^{v\sin\a}\cos(v\cos\a)\frac{\p}{\p G}+
D_2e^{v\sin\a}\sin(v\cos\a)\frac{\p}{\p G},
\end{eqnarray*}
where $A,B\in \C$, and $C, D_1, D_2 \in \R$ are constants\footnote{Note that 
the
corresponding algebra of Lie point symmetries for the heavenly
equation (\ref{pleb1}) 
is infinite-dimensional \cite{BW89}. In order to obtain a finite-dimensional algebra one needs
to factorize it by the infinite-dimensional gauge algebra 
corresponding to the freedom in the definition of $\Om$.
In our case the gauge freedom in $\Om$ was already used to
find the canonical form of the Killing vector. There is no residual
gauge freedom in $F$.}. 
Real generators are
\begin{eqnarray}
X_1&=&\p_w+\p_{\wt},\;\;X_2=i(\p_w-\p_{\wt}),\;\;X_3=i(w\p_w-\wt\p_{\wt})\\
X_4&=&\p_v,\;\;X_5=w\p_w+\wt\p_{\wt}+G\p_G,\nonumber\\
X_6&=&e^{v\sin\a}\sin{(v\cos\a)}\p_G,\;\;
X_7=e^{v\sin\a}\cos{(v\cos\a)}\p_G.\nonumber
\end{eqnarray}
The commutation relations between these vector fields are given by the
following table, the entry in row $i$ and column
$j$ representing $[X_i, X_j]$.
{\small
\begin{center}
\begin{tabular}{p{1cm}|lllllll}
\multicolumn{7}{c}{}\\ 
&$X_1$&$X_2$&$X_3$&$X_4$&$X_5$&$X_6$&$X_7$\\
\hline\
$X_1$&0&0&$X_2$&0&$X_1$&0&0\\
$X_2$&0&0&$-X_1$&0&$X_2$&0&0\\
$X_3$&$-X_2$&$X_1$&0&0&0&0&0\\
$X_4$&0&0&0&0&0&$\sin\a X_6+\cos\a X_7$&$\sin\a X_7-\cos\a X_6$\\
$X_5$&$-X_1$&$-X_2$&0&0&0&$-X_6$&$-X_7$\\
$X_6$&0&0&0&$-\sin\a X_6-\cos\a X_7$&$X_6$&0&0\\
$X_7$&0&0&0&$-\sin\a X_7+\cos\a X_6$&$X_7$&0&0\\
\end{tabular}
\end{center}
\small}
This list of symmetries may seem disappointingly small (as equation
(\ref{alterform}) is an integrable PDE). Further symmetry properties
reflecting
the existence of infinitely many conservation laws will require the
recursive
procedure of constructing `hidden symmetries'. This will be developed
in Section \ref{hidden_symmetries}.
\subsection{Group invariant solutions}
\label{group_solutions} 
 We can  simplify equation (\ref{alterform}) by looking at group
invariant solutions. The finite transformation generated by $X_7$ does not
change the metric. The one by $X_5$ rescales it by a constant factor.
All transformations are conformal Killing vectors for $h$.
\begin{itemize}
\item $X_3=i(w\p_w-\wt\p_{\wt})$ and the corresponding solutions
depend on $(v, R:=\mbox{ln}(w\wt))$. 
This will lead to a new 2D integrable system (\ref{Mequation}). 
Multiplying  (\ref{alterform})  by $e^R$ yields 
\[
(G+G_{vv}-2\sin{\a} G_v)G_{RR}-(e^{i\a}
G_R-iG_{vR})(e^{-i\a}G_R+iG_{vR})=4e^R.
\]
The ideal (\ref{ideal}) reduces to 
\begin{eqnarray*}
0&=&i\d Q\wedge \d J+e^{-i\a}(J\d Q\wedge \d v-Q\d J\wedge \d v)-4\d
(e^R) \wedge \d v,\\
0&=&\d Q\wedge \d v-e^{i\a} J\d R\wedge \d v-i\d J\wedge \d R,
\end{eqnarray*} 
where $J=G_R,\;\;Q=(e^{i\a} G-iG_v)$.
Eliminate $Q$ and use $(J,v)$ as coordinates  to
obtain{\footnote{
With the definition $\xi:=\ln{J}, M:=M(v, \xi)=R-2\xi$ we have
\be
\label{Mequation}
M_{vv}+2M_{v\xi}\sin\a+M_{\xi\xi}+ 4e^M(M_{\xi\xi}+{M_\xi}^2+3M_{\xi}+2)=0.
\ee }} an equation for $R(J,v)$
\be
4(e^R)_{JJ}+R_{vv}+2(JR_J)_v\sin\a+J(JR_J)_J=0.
\ee
Putting $R(J, v)=f(J)+g(v)$ yields (for constant $\a_i$) 
\[
R(J, v)=\a_1 v^2+\a_2 v+\a_1\mbox{arctanh}{\sqrt{4J^{-2}+1}}+\a_3.
\]

A simple solution to (\ref{alterform}) is 
\be
\label{rozw2}
G=e^{v\sin\a}\frac{w\wt}{b}+4e^{-v\sin\a}\frac{b}{1+3\sin^2\a}.
\ee
It has
\[
\Om(w, z, \wt, \zt)=(z\zt)^{(\cos^2\a)/2}
\Big(\frac{\zt}{z}\Big)^{(i\sin\a\cos\a)/2}  \frac{w\wt}{b}+
(z\zt)^{1+(\cos\a)/2}\frac{4b}{1+3\sin^2\a}.
\]
Calculation of curvature components shows it describes a  flat metric
on $\R^4$. Therefore the corresponding EW metric belongs to a class 
described in \cite{PT93}.
\item $X_4=\p_v$. Equation (\ref{alterform}) reduces to
$
GG_{w\wt}-G_wG_{\wt}=4.
$
Define $\Psi(w, \wt)$ by $e^{\Psi}=G$. The EW structure is 
(after rescaling by $16e^{-2\Psi}$) given by
\begin{eqnarray*}
h&=&16e^{-2\Psi}\d w\d \wt+(\d v-ie^{-i\a}\Psi_{w}\d w+
ie^{i\a}\Psi_{\wt}\d \wt)^2,\\ 
\nu&=&-2\sin{\a}\d v+
(2\sin^2{\a}+i\sin{2\a})\Psi_{w}\d w+
(2\sin^2{\a}-i\sin{2\a})\Psi_{\wt}\d \wt,\\
& &\mbox{where}\qquad \Psi_{w\wt}=4e^{-2\Psi}\qquad\mbox{(Liouville equation)}.
\end{eqnarray*}
The general solution to the Liouville equation is
\[
e^{\Psi}=\frac{i(P-\ov{P})}{4\sqrt{P_w\ov{P}_{\wt}}},
\] 
where $P(w)$ is an arbitrary holomorphic function.
With no loss of generality we can take
\be
\label{rozw1}
\Psi=\log{(4b+\frac{w\wt}{b})},\qquad b=const.
\ee
Define new coordinates $(\phi, \theta, \psi)$ by
\[
w=2b\tan{(\theta/2)}e^{i\phi},\qquad \d v=\cos{\a}(\d\psi-\d \phi)+
(\sin{\a})\tan{(\theta/2)}\d\theta
\]
to obtain
\be
\label{TBerger}
h=\d \theta^2+\sin^2{\theta}\d\phi^2+\cos^2{\a}(\d
\psi-\cos{\theta}\d\phi)^2,
\qquad \nu=-\sin{2\a}(\d
\psi-\cos{\theta}\d\phi)
\ee
which is the EW structure on the Berger sphere. 
Calculating the curvature components shows
that the corresponding hyper-K{\"a}hler metric is flat.
The transformation of solution (\ref{rozw1}) corresponding 
to Lie point symmetries
\[
G(w,\wt)\longrightarrow \widetilde{G}(w,\wt, v)=Be^{-v\sin\a}(G(w, \wt)+g(v)),
\]
where
\[
g(v)=-4b+\frac{b e^{2v\sin\a}} {B\cos^2\a}+Ce^{iv\cos\a}+\ov{C}e^{-iv\cos\a}
\]
gives a new solution. In particular (\ref{rozw2}) can be obtained in
this way. Therefore the metric corresponding to (\ref{rozw2})
also describes a Berger sphere.
If $\a=0$ then (\ref{rozw1}) and (\ref{rozw2})
coincide and give the standard metric on $S^3$.
\item $X_2=i(\p_w-\p_{\wt})$ (or $X_1$). This reduction  leads to 
a linear equation.
Put  $w+\wt=f$ to obtain
\[
(G+G_{vv}-2\sin{\a} G_v)G_{ff}-(e^{i\a}
G_R-iG_{vf})(e^{-i\a}G_f+iG_{vf})=4.
\]
With the definition $J:=G_f,\;\;Q:=(e^{i\a} G-iG_v)$
this yields
\begin{eqnarray*}
0&=&i\d Q\wedge \d J+e^{-i\a}(J\d Q\wedge \d v-Q\d J\wedge \d v)-4\d f
\wedge \d v,\\
0&=&\d Q\wedge \d v-e^{i\a} J\d f\wedge \d v-i\d J\wedge \d f.
\end{eqnarray*}
Now eliminate $Q$ and use $(v, \xi=\ln{J})$ as coordinates  to
obtain a linear equation for $f(\xi,v)$
\be
4e^{-2\xi}(f_{\xi \xi}-f_{\xi})+f_{vv}+2\sin\a f_{\xi v}+f_{\xi \xi}=0.
\ee
\end{itemize}
\section{Hidden symmetries}
\label{hidden_symmetries}
In this section we shall find a recursion procedure
for generating `hidden symmetries' of (\ref{euclid}).
We start with discussing the general conformally invariant wave
equation in Einstein--Weyl background.

A tensor object $T$ which transforms as 
\[
T\longrightarrow \phi^m T\qquad\mbox{when}\qquad 
h_{ij}\longrightarrow \phi^2 h_{ij} 
\]
is said to be conformally invariant of weight $m$.
Let $\beta$ be a $p$-form of weight $m$. The covariant
derivative
\[
D\beta:=\d\beta-\frac{m}{2}\nu\wedge\beta
\]
is a well defined $p+1$ form of weight $m$. Its Hodge dual, $\ast_h
D\beta$,
is a $(2-p)$-form of weight $m+1-p$. Therefore we  can write the
weighted Weyl
wave operator which takes $p$-forms of weight $m$ to $(3-p)$-forms
of weight $m+1-p$
\[
D\ast_hD=\Big(\d-\frac{m+1-p}{2}\nu\wedge\Big)\ast_h
\Big(\d-\frac{m}{2}\nu\wedge\Big).
\]
Consider the case $p=0$. Let $\phi$ be a function 
of weight $m$. The most general wave equation is
\[
D\ast_hD\phi=kW\phi\;{vol}_h
\]
where $k$ is some constant. The RHS  has weight $m+1$ 
so the whole expression is conformally invariant.
Adopting the index notation we obtain
\be
\label{confequation}
\nabla^i\nabla_i\phi-\Big(m+\frac{1}{2}\Big)\nu^i\nabla_i\phi
+\frac{1}{4}\Big(m(m+1)\nu^i\nu_i-2m\nabla^i\nu_i\Big)\phi=
k\Big(R+2 \nabla^i\nu_i-\frac{1}{2}\nu_i\nu^i\Big)\phi.
\ee 
At this stage one can make some choices concerning the values of 
$m$ and $k$. One
can also fix the gauge freedom. In \cite{CTV96} it was assumed that
$k=0, m=-1$ and $\nabla_i\nu^i=0$ (the  Gauduchon gauge)
which led to the derivative of the generalised
monopole equation (\ref{EWmonopole}):
\[
\nabla^i\nabla_i\phi+\frac{1}{2}\nu^i\nabla_i\phi=0.
\]
Another possibility is to set $m=-(1/2), k=1/8$. 
With this choice equation  (\ref{confequation}) simplifies to
\[
\nabla^i\nabla_i\phi=\frac{1}{8}R\phi,
\]
which is the well known conformally invariant wave equation in the
3D Riemannian geometry. Note  that the gauge freedom was not fixed
to derive the last equation. All we did was to get rid of
the `non-Riemannian' data. 
\subsection{The recursion procedure}
Let $\delta F$ be a linearised solution to  (\ref{EWred})
(i.e. $F+\delta F$ satisfies  (\ref{EWred}) up to the linear terms
in $\delta F$). Then
\begin{eqnarray}
\label{EWlinear}
&&\Big([(\eta F_{w}-F_{uw})\frac{\p^2}{\p u\p {\tw}}
-(\eta F_{\tw}+F_{u\tw})\frac{\p^2}{\p u\p w} 
-(\eta^2F-F_{uu})\frac{\p^2}{\p \tw\p w}
+F_{w\tw}\frac{\p^2}{\p u^2}]\nonumber\\
&&+\eta[(\eta F_{w}-F_{uw})\frac{\p}{\p {\tw}}
+\eta(\eta F_{\tw}+F_{u\tw})\frac{\p}{\p w}] 
\Big)\delta F
=F_{w\tw}\delta F.
\end{eqnarray}
This equation can be viewed more geometrically:
let $\square_{\Om}$ denote the wave operator on an 
ASDV curved background given by
$\Om$, let $\dom$ be the linearised solution to the first
heavenly equation
and let ${\cal W}_{\Om}$ be the kernel of  $\square_{\Om}$.
It is straightforward to check \cite{DM96} that $\dom \in {\cal W}_{\Om}$.
Indeed, put $\p:=\d w\otimes\p_w+\d z\otimes\p_z$, $\pt:=\d
\tw\otimes\p_{\tw}+
\d {\tz}\otimes\p_{\tz}$
and rewrite (\ref{pleb1}) as $(\p \pt (\Om +\dom ))^2 =\nu$. For 
the linearised solution we have
\[
0=(\p \pt \Om\wedge \p \pt ) \dom = \d( \p \pt \Om\wedge (\p -\pt)\dom )=
\d\ast_g \d\dom=\square_{\Om}\dom.
\]
Now impose the additional constrain ${\cal L}_K\dom=\eta \dom$. This
implies $\dom=e^{\eta t}\delta F$. This yields
\begin{eqnarray*}
0&=&\d\ast_g\d(e^{\eta t}\delta F)\\
&=&e^{\eta t}
((\eta^2(\d t\wedge\ast_g\d t)+\eta \d
\ast_g\d t)\delta F + \eta\d t\wedge \ast_g\d \delta F+
\eta\d \delta F\wedge\ast_g\d t +\d \ast_g\d \delta F).
\end{eqnarray*}
But $\d \ast_g\d t=\square_{\Om} t=0$ and 
$\d t\wedge\ast_g\d t=|\d t|^2\nu_g$,
therefore (\ref{EWlinear}) is equivalent to
\[
\square_{\Om} \delta F+\eta^2|\d t|^2\delta F=0.
\]
There should exist a choice
of $m$ and $k$ which, in the appropriate gauge, reduces equation
(\ref{confequation}) down to (\ref{EWlinear}).

 Let ${\cal W}_{F}$ be the space of solutions to (\ref{EWlinear})
around a given solution $F$. We shall construct a map
$R:{\cal W}_{F}\longrightarrow {\cal W}_{F}$. 
Let us
start from the recursion operator for the heavenly equation
\cite{DM96}. Let $\phi\in  {\cal W}_{\Om}$.
Define a recursion operator $R:{\cal W}_{\Om}\longrightarrow {\cal
W}_{\Om}$
by
\be
\label{def}
o^{A'}\nabla_{AA'}R\phi=\iota^{A'}\nabla_{AA'}\phi,\;\;\;\;
o^{A'}=(1,0),\;\iota^{A'}=(0, 1),\;\;\; A=0,1,
\ee
where, in coordinates $(w, \tw, t, u)$
\begin{eqnarray*}
\nabla_{00'}&=& \frac{m}{2}e^{\r t+\tilde{m}u}\Big(F_{w\tw}
\Big(\frac{\p}{\p t}-\frac{\p}{\p u}\Big)-(\eta
F_w-F_{uw})\frac{\p}{\p{\tw}}\Big),\\
\nabla_{10'}&=&\frac{1}{4}e^{-\eta t-2\r u}\Big((\eta
F_{\tw}+F_{u\tw})
\Big(\frac{\p}{\p t}-\frac{\p}{\p u}\Big)   
-(\eta^2 F-F_{uu})\frac{\p}{\p{\tw}}\Big),\\
\nabla_{01'}&=&\frac{\p}{\p w},\\
\nabla_{11'}&=&\frac{1}{2m}e^{-m(t+u)}\Big(\frac{\p}{\p t}+\frac{\p}{\p u}\Big).
\end{eqnarray*}
To construct a reduced recursion operator we 
should be able to  Lie derive (\ref{def}) along $K$. In order to
do so we introduce an invariant spin frame 
\[
\hat{o}^{A'}:=e^{-(1/2)\r t}o^{A'},\;\;
\hat{\iota}^{A'}:=e^{(1/2)\r t}\iota^{A'},
\]
in which $\lt=(\pi_{A'}\hat{o}^{A'})/(\pi_{A'}\hat{\iota}^{A'})$.
Note that now  $\Gamma_{A'B'}\neq 0$.
Recursion relations are
\[
e^{-\r t}\nabla_{A0'}(e^{\eta t}R\delta F)=\nabla_{A1'}e^{\eta t}\delta 
F.
\]
This yields the following result
\begin{prop}
The map $R:{\cal W}_{F}\longrightarrow {\cal W}_{F}$ defined by
\begin{eqnarray}
\label{EWrecur}
me^{\tilde{m}u}(F_{w\tw}(\eta-\p_u)-(\eta F_w-F_{uw})\p_{\tw})
R\delta F&=&2\p_w\delta F\\
\tilde{m}e^{\tilde{m}u}((\eta F_{\tw}+F_{u\tw})(\eta-\p_u)
-(\eta^2F-F_{uu})\p_{\tw})R\delta F&=&2(\eta+\p_u)\delta F.\nonumber
\end{eqnarray}
generates new elements of  ${\cal W}_{F}$
from the old ones.
\end{prop}
By cross differentiating we verify that two equations in (\ref{EWrecur})
are consistent as a consequence of (\ref{EWred}).
 
 We start the recursion from two solutions $(e^{-\eta u},
\frac{2\tilde{m}}{m+\eta}e^{mu})$ to (\ref{EWlinear}).
Equations (\ref{EWrecur}) yield
\[
e^{-\eta u}\longrightarrow -\frac{\eta F+F_u}{2m}\longrightarrow...,\;\;\;\;
\frac{2\tilde{m}}{m+\eta}e^{mu}\longrightarrow F_w\longrightarrow...\;.
\]
Suppose that $F=F(u, w, \tw, {\bf T})$ depends on three local coordinates on
a complex EW space and a sequence of parameters ${\bf T}=(T_2, T_3,
... )$. Put 
\[
\frac{\p F}{\p T_n}:=R^n\Big(\frac{2\tilde{m}}{m+\eta}e^{mu}\Big),
\]
so that $T_1=w$. The recursion relations $R(\p_{T_n}F)=\p_{T_{n+1}}F$
form an over-determined  system of equations which involve
arbitrarily many independent variables, but initial data can be
specified freely only on a two dimensional surface.

\section{Alternative formulations}
\label{alternative_formulations}
Here we shall give an alternative formulation 
of equation (\ref{EWred}).
Define functions $(V, S, \widetilde{S})$ by
\[
4V:=\eta^2F-F_{uu},\;\;\;\;\;\; 2S:=\eta F_{w}-F_{uw},\;\;
2\widetilde{S}:=\eta F_{\tw}+F_{u\tw},
\]
so equation (\ref{EWred}) takes the form
\be
\label{todform}
V=\frac{(-e^{2\r u}+S\widetilde{S})\eta}{S_{\tw}+\widetilde{S}_{w}},
\;\;\;\ S_u+\eta S=2V_w,\;\;\;-\widetilde{S}_u+\eta \widetilde{S}=2V_{\tw}.
\ee
The hyper-K{\"a}hler metric is
\begin{eqnarray*}
\d s^2&=&e^{\eta t}(V(\d t^2-\d u^2)+V^{-1}(S\widetilde{S}-e^{2\r u})\d
w\d \tw+S(\d t-\d u)\d w+\widetilde{S}(\d t+\d u)\d \tw\\
&=&e^{\eta t}(V^{-1}h+V(\d t+\om)^2)
\end{eqnarray*}
where
\[
h:=-e^{2\r u}\d w\d \tw -\Big(V\d u+\frac{S\d w -\tilde{S} \d \tw}{2}\Big)^2
\;\;\;\ \om:=\frac{S\d w +\tilde{S} \d \tw}{2V},
\]
and the EW one-form corresponding to $h$ is
\[
\nu
=4\rho\d u+\frac{(\eta+2\rho)S\d w+(\eta-2\rho)\tilde{S} \d \tw}{V}.
\]
Euclidean reality conditions force $\ov{S}=-\widetilde{S}$ and $V$
real. On the $++--$ slice we have $\ov{S}=\widetilde{S}$, or alternatively
(on a different real slice) 
functions $V, S, \widetilde{S}$  real and independent.
The orthonormal frame on the Euclidean slice is
\begin{eqnarray}
\label{euclidtetrad}
e^1&=&\frac{1}{2}(e^{imv}\d w+e^{-i\ov{m}v}\d \wt),\;\;\;\;\;
\nabla_1=e^{-imv}\p_w+e^{i\ov{m}v}\p_{\wt}
+i\frac{Se^{-imv} -\overline{S}e^{i\ov{m}v} }{2V}\p_v\\
e^2&=&\frac{i}{2}(e^{-i\ov{m}v}\d \wt-e^{imv}\d w ),\;\;\;\;\;
\nabla_2=i(e^{-imv}\p_w-e^{i\ov{m}v}\p_{\wt})
-\frac{Se^{-imv} +\overline{S}e^{i\ov{m}v} }{2V}\p_v\nonumber\\
e^3&=&V\d v-i\frac{S\d w -\overline{S} \d \wt}{2},\;\;\;\;\;\;\;\;
\nabla_3=\frac{1}{V}\p_v.\nonumber\\
\end{eqnarray}
The EW one form is
\[
\nu=2\om\cos{\a}-\frac{4\sin{\a}}{V}e^3=\frac{\cos{\a}(S\d w+\overline{S} \d
\wt)-2i\sin{\a}(\overline{S} \d \wt-S\d w)}{V}-4\sin{\a}\d v.
\]
Equations (\ref{todform}) can be rewritten in a compact form
\be
\label{af1}
\d e^3=\om\wedge e^3\cos\a+\frac{\cos\a}{V}e^1\wedge e^2,
\ee
\be
\label{af2}
\d(e^1+ie^2)=e^{i\a}\om\wedge(e^1+ie^2)+\frac{ie^{i\a}}{V}e^3
\wedge(e^1+ie^2)
\ee
(the last relation is an identity).
In fact the converse is true:
\begin{prop}
Let $(e^1, e^2, e^3)$ be real one-forms which satisfy 
{\em (\ref{af1},\ref{af2})}
for some real one-form $\om$, function $V$ and constant $\a$. Then there
exist local coordinates $w\in \C, v\in\R$ and a complex function
$S(w, \wt, v)$ such that $(e^1, e^2, e^3)$ are of the form
{\em(\ref{euclidtetrad})} and the Euclidean version of {\em(\ref{todform})}
is satisfied.
\end{prop}
{\bf Proof.}
Equation (\ref{af2}) and the Frobenuis theorem imply that 
$e^1+ie^2=e^{im\chi}\d w$  (where $m=e^{i\a}$)
for some complex functions $\chi$ and $w$, which
therefore satisfie
\[
\d \chi =\om+\frac{i}{V}e^3-S\d w
\]
for some $S$. Put $\chi= v+iY$ (for $v,Y\in\R$) so that
\[
\d v=\frac{1}{V}e^3+\frac{i}{2}(S\d w-\ov{S}\d \wt),\qquad
\d Y=\frac{1}{2}(S\d w+\ov{S}\d \wt)-\om.
\]
Now we use the conformal freedom of (\ref{af1},\ref{af2}) and rescale
\[
\hat{e}^1+i\hat{e}^2=\Phi e^{i\theta}(e^1+ie^2),\qquad
\hat{e}^3=\Phi e^3,\qquad\hat{V}=\Phi V,\]
\[
\hat{\om}=\om+ (\cos{\a}\Phi)^{-1}\d \Phi,
\qquad\d\theta=-\tan{\a}\Phi^{-1}\d \Phi,
\]
so that we can put $Y=0$, and (\ref{af2}) is solved.
Now  
\[
\om=\frac{S\d w+\ov{S}\d\wt}{2V},\qquad
e^3=V\d v-i\frac{S\d w -\overline{S} \d \wt}{2},
\]
and the equation (\ref{af1}) gives (\ref{todform}).
\koniec
Recall that a geodesic congruence $\Gamma$ in a region 
$U\subset{\cal W}$ is a set of geodesics, one through each point of $U$.
Let $W^i$ be a generator of $\Gamma$ (a vector field 
tangent to $\Gamma$).
Then the geodesic condition is $W^jD_jW^i\sim W^i$.
The formula (\ref{af1}) implies that
$e^3$ generates a shear--free geodesic congruence, 
with  twist and divergence
given by:
\[
\mbox{twist}=*_h(e^3\wedge\d e^3)=\frac{\cos{\a}}{V},\qquad
\mbox{divergence}=*_h\d*_he^3=\frac{6\sin{\a}}{V}.
\]
They are both solutions 
of the generalized monopole equation (\ref{EWmonopole}).
Conversely, it follows from \cite{CP99} that if the twist and the 
divergence 
of a shear--free geodesic congruence on an EW space are proportional,
then this EW space arises as a reduction of a hyper-K\"ahler 
metric\footnote{We are grateful to David Calderbank for informing 
us of the results
in \cite{CP99}.}.
Therefore solutions 
to (\ref{alterform}) (or equivalently the Euclidean version of
(\ref{todform})) are completely characterized by the existence of 
a shear--free geodesic congruence of the above type.
It should however be stressed that, given an EW structure, there is no
a priori way of telling if this special shear--free geodesic congruence exists.
It would be interesting to find a local obstruction 
to the existence of such congruence.


\section{Acknowledgements}
We are grateful to David Calderbank and Lionel Mason for useful discussions.
This work was finished during the workshop
{\em Spaces of geodesics and complex methods in general relativity and 
geometry } held in the summer of 1999 at the Erwin Schr{\"o}dinger 
Institute in Vienna.  We wish to thank ESI  
for the hospitality and for financial
assistance.


\end{document}